\documentclass[12pt]{amsart}
\usepackage{amsmath}
\usepackage{amssymb}
\usepackage{latexsym}

\newcommand{\Zint}{\mathbb {Z}}    
     
\newcommand{\Rea}{\mathbb {R}}      

\newcommand{\halmos}{\rule{5pt}{5pt}}

\numberwithin{equation}{section}
\newtheorem{dfn}{\bf Definition}
\newtheorem{prop}{\bf Proposition}[section]

\newtheorem{thm}[prop]{\bf Theorem}

\newtheorem{exa}{\bf Example}

\begin{document}

\title[Finite-gap potential]
{Finite-gap potential, Heun's differential equation and WKB analysis}
\author{Kouichi Takemura}
\address{Department of Mathematical Sciences, Yokohama City University, 22-2 Seto, Kanazawa-ku, Yokohama 236-0027, Japan.}
\email{takemura@yokohama-cu.ac.jp}


\begin{abstract}
We review several results on the finite-gap potential and Heun's differential equation, and we discuss relationships among the finite-gap potential, the WKB analysis and Heun's differential equation.
\end{abstract}

\maketitle

\section{Introduction}

Heun's equation (Heun's differential equation) is a linear differential equation of second order given by
\begin{align}
& \frac{d^2y}{dz^2} \! + \left( \frac{\gamma}{z}+\frac{\delta }{z-1}+\frac{\epsilon}{z-t}\right) \frac{dy}{dz} +\frac{\alpha \beta z -q}{z(z-1)(z-t)} y =0,
\label{Heun}
\end{align}
with the condition $\gamma +\delta +\epsilon =\alpha +\beta +1$ \cite{Ron}.
It has four singularities $\{ 0, 1, t, \infty \}$ and they are all regular.
Heun's equation is known to be a standard form of the second-order Fuchsian differential equation with four singularities.
The parameter $q$ is not determined by the local monodromy, and is called an accessory parameter.
Heun's differential equation frequently appears in Physics, i.e. black hole (general relativity, Kerr's solution), crystalline materials \cite{SL}, fluid dynamics \cite{CH}, quantum mechanics (Inozemtsev model \cite{Ino}) and so on.

A standard form of the second-order Fuchsian differential equation with three singularities is given by the hypergeometric differential equation
\begin{equation}
z(1-z) \frac{d^2y}{dz^2} \! + \left( \gamma - (\alpha + \beta +1)z \right) \frac{dy}{dz} -\alpha \beta y=0.
\end{equation}
It has three regular singularities at $\{ 0, 1, \infty \}$.
Global properties of solutions and the monodromy are known for the hypergeometric differential equation.
In particular an integral representation of a solution is given by
\begin{align}
& F(\alpha ,\beta ;\gamma ;z) = \frac{\Gamma (\gamma )}{\Gamma (\alpha )\Gamma (\gamma -\alpha )}\int _0 ^1 s^{\alpha -1}(1-s)^{\gamma - \alpha -1}(1-sz)^{- \beta } ds.
\end{align}
The connection matrix among local solutions at $z=0$ and ones at $z=1$ is written in terms of the gamma function.
The condition for existence of non-zero function holomorphic at $z=0$ and $z=1$ is known and the holomorphic function is given by the Jacobi polynomial.

For investigating global properties for Heun's differential equation, we can apply the method of finite-gap integration.
For this purpose, we recall an elliptic representation of Heun's differential equation.
Let $\wp (x)$ be the Weierstrass doubly-periodic function with periods $(2\omega _1 ,2\omega _3)$. 
Set $\omega _0=0, \omega _2=-\omega _1 -\omega _3$, $e_i=\wp(\omega_i) \; (i=1,2,3)$,
\begin{align}
& z=\frac{\wp (x) -e_1}{e_2-e_1}, \quad  t=\frac{e_3-e_1}{e_2-e_1}, \quad f(x)= y z^{\frac{-l_1}{2}}(z-1)^{\frac{-l_2}{2}}(z-t)^{\frac{-l_3}{2}},
\end{align}
then Heun's differential equation (Eq.(\ref{Heun})) is transformed to
\begin{equation}
\left(-\frac{d^2}{dx^2} + \sum_{i=0}^3 l_i(l_i+1)\wp (x+\omega_i)-E\right)f(x)=0,
\label{InoEF0}
\end{equation}
where 
\begin{align}
& l_0= \beta -\alpha -1/2,\quad l_1= -\gamma +1/2, \quad l_2=-\delta +1/2, \quad l_3=-\epsilon +1/2, \\
& E=(e_2-e_1)(-4q+(-(\alpha -\beta)^2 +2\gamma ^2+6\gamma \epsilon +2\epsilon ^2 -4\gamma -4\epsilon  -\delta ^2 +2\delta +1)/3 \nonumber \\
& \quad \quad +(-(\alpha -\beta ) ^2 +2\gamma ^2+6\gamma \delta +2\delta^2 -4\gamma -4\delta -\epsilon ^2+2\epsilon +1)t/3). \nonumber
\end{align}
If $l_1=l_2=l_3=0$ ($\gamma =\delta = \epsilon =1/2$), then the differential equation is called Lam\'e's equation.

It is known that, if $l_0, l_1, l_2, l_3 \in \Zint$, then the function $\sum_{i=0}^3 l_i(l_i+1)\wp (x+\omega_i)$ is an algebro-geometric finite-gap potential, and is called the Treibich-Verdier potential.
In section \ref{sec:fingap}, we recall the definitions on the finite-gap potential and review the properties and the examples including the Treibich-Verdier potential.
The monodromy of Heun's differential equation is investigated in connection with the finite-gap property, and we will obtain formulae related with the finite-gap property in section \ref{sec:monod}.
In section \ref{sec:WKB}, we discuss relationships among the finite-gap potential, the WKB analysis and Heun's differential equation.
In particular, we provide another approach for results by Borcea and Shapiro \cite{BS} on root asymptotics of spectral polynomials for the Lam\'e operator.

\section{Finite-gap potential} \label{sec:fingap}

We recall definitions of the finite-gap potential and the algebro-geometric finite-gap potential.

\begin{dfn}
Let $q(x)$ be a periodic, smooth, real function, $H$ be the operator $-d^2/dx^2+q(x)$, and the set $\sigma _b(H)$ be defined as follows:
\begin{equation}
E \in \sigma _b(H) \; \Leftrightarrow \mbox{ Every solution to }(H-E)f(x)=0 \mbox{ is bounded on }x \in \Rea .
\end{equation}
If the closure of the set $\sigma _b(H)$ can be written as
\begin{equation}
\overline{\sigma _b(H)}= [E_{0},E_{1}]\cup [E_{2},E_{3}]  \cup \dots \cup [E_{2g}, \infty ),
\end{equation}
where $E_0<E_{1}<\cdots <E_{2g}$, i.e., the number of bounded bands is finite, then $q(x)$ is called the finite-gap ($g$-gap) potential.
\end{dfn}

\begin{exa}
We consider the case $q(x)=0$, i.e., $H=-d^2/dx^2$.
Fix the eigenvalue $E$. Then $H f(x)=E f(x)$ is equivalent to $f''(x)+E f(x)=0$.
We solve the differential equation for dividing into three cases.
If $E<0$, then we write $E=-\lambda ^2$ and the solutions are $f(x)=Ae^{\lambda x} +B e^{-\lambda x} $ for constants $A,B$, which are unbounded on $\Rea$ for $(A,B)\neq (0,0)$.
If $E=0$, then the solutions are $f(x)=A+B x $, which are unbounded on $\Rea$ for $B\neq 0$.
If $E>0$, then write $E=\lambda ^2$ and the solutions are $f(x)=A\cos \lambda x +B \sin \lambda x $, which are bouned on $\Rea$.
Hence we have $\sigma _b(H)= (0,\infty)$ and $\overline{\sigma _b(H)}= [0,\infty)$.
Therefor the potential $q(x)=0$ is finite-gap ($0$-gap).
\end{exa}

Assume that the potential $q(x)$ is real, smooth, periodic with a period $T$.
Let $f_1(x,E)$, $f_2 (x,E)$ be a basis of solutions to $\left ( -d^2/dx^2+q(x) -E\right) f(x)=0$.
Then $f_1(x+T,E)$ and $f_2 (x+T,E)$ are also solutions and written as
\begin{align}
& (f_1(x+T,E) \; f_2(x+T,E)) = (f_1(x,E) \; f_2(x,E))M,
\end{align}
where $M$ is a $2\times 2$ matrix with constant elements. It is known that $\det M=1$.
Let $t^2-($tr$ M) t+1=0$ be the characteristic polynomial of the monodromy matrix $M$.
If $|$tr$M |<2$ (resp. $|$tr$M |>2$), then $E \in \sigma _b(H)$ (bounded) (resp. $E \not \in \sigma _b(H)$ (unbounded)).
If tr$M =2$ (resp. tr$M =-2$), then there exists a non-zero periodic solution $f(x+T)=f(x)$ (resp. an antiperiodic solution $f(x+T)=-f(x)$).
Hence the monodromy caused by the shift of the period $(x \mapsto x+T)$ implies the boundedness or the unboundedness of the solutions to the differential equation.

The definition of algebro-geometric finite-gap potential is described as follows:
\begin{dfn}
If there exists an odd-order differential operator 
\begin{equation}
A= \left( \frac{d}{dx} \right)^{2g+1} +  \sum_{j=0}^{2g-1}\! b_j(x) \left( \frac{d}{dx} \right)^{2g-1-j} 
\end{equation}
such that 
\begin{equation}
\left[ A, -\frac{d^2}{dx^2}+q(x)\right] =0,
\end{equation}
then $q(x)$ is called the algebro-geometric finite-gap potential.
\end{dfn}
Note that the equation  $[A, -d^2/dx^2+q(x)]=0$ is equivalent to the function $q(x)$ being a solution to a stationary higher-order KdV equation (see \cite{DMN}).
It was established in the 1970s that, under the assumption that $q(x)$ is a periodic, smooth, real function, the potential $q(x)$ is finite-gap if and only if $q(x)$ is algebro-geometric finite-gap.

We now present examples of the (algebro-geometric) finite-gap potentials.
Ince \cite{I} established in 1940 that if $n \in \Zint _{\geq 1}$, $\omega _1 \in \Rea \setminus \{ 0 \}$ and $\omega _3 \in \sqrt{-1} \Rea \setminus \{ 0 \}$, then the potential of Lam\'e's operator 
\begin{equation}
-\frac{d^2}{dx^2}+n(n+1)\wp (x+\omega _3),
\end{equation}
is finite-gap.
In the late 1980s, Treibich and Verdier \cite{TV} found that the method of finite-gap integration is applicable the elliptic representation of Heun's equation for the case $l_0 , l_1 , l_2 , l_3 \in \Zint $.
Namely, they showed that the potential in Eq.(\ref{InoEF0}) is an algebro-geometric finite-gap potential if $l_0, l_1, l_2, l_3 \in \Zint $.
Therefore the potential $ \sum_{i=0}^3 l_i(l_i+1)\wp (x+\omega_i)$ is called the Treibich-Verdier potential.
Treibich and Verdier developed the theory of elliptic soliton, Jacobi variety and tangential covering, and obtained the results.
Subsequently several others \cite{GW,Smi,Tak1,Tak2,Tak3,Tak4,Tak5} have produced more precise statements and concerned results on this subject.

We obtained further examples of algebro-geometric finite-gap potentials in \cite{TakF}.
\begin{thm} (\cite{TakF}) \label{thm:fingap} 
If $M, l_0 , l_1, l_2 , l_3 \in \Zint _{\geq 0}$, $\delta _j \not \equiv \omega _i$ mod $2\omega_1 \Zint \oplus 2\omega_3 \Zint$ $(0\leq i\leq 3, \; 1\leq j\leq M)$, $\delta _j \pm \delta _{j'} \not \equiv 0$  mod $2\omega_1 \Zint \oplus 2\omega_3 \Zint$ $(1\leq j< j' \leq M)$ and $\delta _1 ,\dots, \delta _M$ satisfy the equations
\begin{equation}
2\sum _{j' \neq j} (\wp ' (\delta _j -\delta _{j'} ) + \wp ' (\delta _j +\delta _{j'} ) ) +\sum _{i=0}^3 (l_i +1/2)^2 \wp ' (\delta _j +\omega _i ) =0,
\label{eq:cond} 
\end{equation}
$(j=1,\dots ,M)$, then the potential
\begin{align}
v(x) = & \sum_{i=0}^3 l_i(l_i+1) \wp (x+\omega_i) +2\sum_{i'=1}^M (\wp (x-\delta _{i'}) + \wp (x+\delta _{i'})),
\end{align}
is algebro-geometric finite-gap.
\end{thm}
Note that Eq.(\ref{eq:cond}) has appropriate solutions for each $M, l_0 , l_1, l_2 , l_3 \in \Zint _{\geq 0}$ (see \cite{TakF}).
If $M=0$, then we recover the result on Heun's equation, and if $M=1$, then we recover the Treibich's result \cite{Tre}.
Gesztesy and Weikard \cite{GW2} developed a theory of Picard's potential, and it would be related to our one.

We introduce a proposition which plays a crutial role of proving Theorem \ref{thm:fingap}.
Observe that a product of two solutions to 
\begin{equation}
\left( -\frac{d^2}{dx^2}+v(x) \right) f(x)=Ef(x),
\label{eq:DEfingap}
\end{equation}
satisfies
\begin{equation}
 \left\{ \frac{d^3}{dx^3}-4\left( v(x) -E \right) \frac{d}{dx} -2\frac{dv(x)}{dx} \right\} \Xi (x,E)=0.
\label{prodDE}
\end{equation}
\begin{prop} \label{prop:prod} (\cite{TakF}).
Under the condition (\ref{eq:cond}), Eq.(\ref{prodDE}) has a unique non-zero doubly-periodic solution $\Xi (x,E)$, which has the expansion
\begin{align}
& \Xi (x,E)=c_0(E)+\sum_{i=0}^3 \sum_{j=0}^{l_i-1} b^{(i)}_j (E)\wp (x+\omega_i)^{l_i-j} \label{Fx1}\\
& \quad \quad \quad + \sum _{i'=1}^M d^{(i')}(E)(\wp (x+\delta _{i'}) +\wp (x- \delta _{i'})), \nonumber
\end{align}
where the coefficients $c_0(E)$, $b^{(i)}_j(E)$ and $d^{(i')}(E)$ are polynomials in $E$, these polynomials do not share any common divisors, and the polynomial $c_0(E)$ is monic. 
\end{prop}
We set $g=\deg_E c_0(E)$. Then the coefficients satisfy $\deg _E b^{(i)}_j(E)<g$ for all $i$ and $j$.
Note that the function $\Xi (x,E)$ is frequently used for describing solutions to Eq.(\ref{eq:DEfingap}) and the monodromy, as we will see in section \ref{sec:monod}.
Write
\begin{equation}
\Xi(x,E) = \sum_{i=0}^{g} a_{g-i}(x) E^i. 
\label{Xiag}
\end{equation}
Then $a_0(x)=1$ and it follows from Eq.(\ref{prodDE}) that
\begin{equation}
a'''_j(x)-4v(x)a'_j(x)-2v'(x)a_j(x)+4a'_{j+1}(x)=0.
\label{a'''v}
\end{equation}
Define the $(2g+1)$st-order differential operator $A$ by
\begin{align}
& A= \sum_{j=0}^{g} \left\{ a_j(x)\frac{d}{dx}-\frac{1}{2} \left( \frac{d}{dx} a_j(x) \right) \right\} H ^{g-j}, \quad  H=-\frac{d^2}{dx^2}+v(x).
\label{Adef}
\end{align}
It follows from Eq.(\ref{a'''v}) that $[A,H]=0$.
Hence the function $v(x)$ is an algebro-geometric finite-gap potential, and we obtain Theorem \ref{thm:fingap}.
Set
\begin{align}
& Q(E)= \Xi (x,E)^2\left( E- v(x)\right) +\frac{1}{2}\Xi (x,E)\frac{d^2\Xi (x,E)}{dx^2}-\frac{1}{4}\left(\frac{d\Xi (x,E)}{dx} \right)^2.
\label{const}
\end{align}
It is shown by differentiating the right-hand side of Eq.(\ref{const}) and applying Eq.(\ref{prodDE}) that $Q(E)$ is independent of $x$, and $Q(E)$ is a monic polynomial in $E$ of degree $2g+1$.
On the operators $H$ and $A$, we have the relation $A^2+Q(H)=0$ (see \cite[Proposition 3.2]{Tak3}).
\begin{exa}
For the case $M=0$, $l_0=2$, $l_1=l_2=l_3=0$, we have
\begin{align}
& \Xi (x,E) = a_0(x) E^2 +a_1(x) E +a_2(x) = E^2 +3\wp (x) E+ 9(\wp (x)^2 - g_2/4),\\
& A = \frac{d}{dx}\left( -\frac{d^2}{dx^2} +6\wp (x) \right) ^2+3\left( \wp (x) \frac{d}{dx} -\frac{1}{2}\wp ' (x) \right)\left( -\frac{d^2}{dx^2} +6\wp (x) \right) \\
& \quad \quad +9\left\{ \left( \wp (x)^2-\frac{g_2}{4} \right) \frac{d}{dx} -\wp (x) \wp '(x) \right\} \nonumber \\
& \quad = \left( \frac{d}{dx} \right)^{5} -15\wp (x) \left( \frac{d}{dx} \right)^{3}  -\frac{45}{2} \wp ' (x) \left( \frac{d}{dx} \right)^{2} -9\left( 5\wp (x)^2-\frac{3}{4}g_2 \right) \frac{d}{dx}, \nonumber \\
& Q(E)=(E^2 -3g_2)(E^3-9g_2E/4-27g_3/4),
\end{align}
where $g_2= -4(e_1e_2+e_2e_3+e_3e_1)$ and $g_3=4e_1e_2e_3$.
\end{exa}

\section{Monodromy related with finite-gap potential} \label{sec:monod}

In this section, we review results on solutions to the Schr\"odinger equation (see Eq.(\ref{eq:DEfingap})) with the algebro-geometric finite-gap potential and its monodromy.
Namely, we have an integral representation for a solution to Eq.(\ref{eq:DEfingap}), a monodromy formula in terms of hyperelliptic integral, an expression of the Hermite-Krichever Ansatz and hyperelliptic-to-ellitic reduction integral formula.
The following propositions are correct for the potentials in Theorem \ref{thm:fingap}, which include the case of Heun's didderential equation.
\begin{prop} (Integral representation for eigenfunctions (\cite{Tak1} for Heun's equation))
Let $\Xi (x,E)$ be the function defined in Proposition \ref{prop:prod} and $Q(E)$ be the polynomial defined in Eq.(\ref{const}). Then
\begin{equation}
\Lambda (x,E)=\sqrt{\Xi (x,E)}\exp \int \frac{\sqrt{-Q(E)}dx}{\Xi (x,E)}
\label{eqn:Lam}
\end{equation}
is a solution to Eq.(\ref{eq:DEfingap}).
\end{prop}
\begin{prop} (Monodromy formula in terms of hyperelliptic integral (\cite{Tak3} for Heun's equation)) \label{thm:conj3} 
Let  $k \in \{ 1,3\}$, $q_k \in \{0,1\}$ and $E_0$ be the eigenvalue such that $\Lambda (x+2\omega _k,E_0)=(-1)^{q_k} \Lambda (x,E_0)$. Then
\begin{align} 
& \Lambda (x+2\omega _k,E)=(-1)^{q_k} \Lambda (x,E) \exp \left( -\frac{1}{2} \int_{E_0}^{E}\frac{ \int_{0+\varepsilon }^{2\omega _k+\varepsilon }\Xi (x,\tilde{E})dx}{\sqrt{-Q(\tilde{E})}} d\tilde{E}\right).
\label{analcontP}
\end{align}
\end{prop}
Write
\begin{align}
\Xi (x,E)= & \: c(E)+\sum_{i=0}^3 \sum_{j=0}^{l_i-1 } a^{(i)}_j (E)\left( \frac{d}{dx} \right) ^{2j} \wp (x+\omega_i) \\
& +\sum _{i'=1}^M d^{(i')}(E)(\wp (x+\delta _{i'}) +\wp (x- \delta _{i'})) , \nonumber \\
a(E)= & \sum _{i=0}^3 a^{(i)} _0 (E)+2 \sum _{i'=1}^M d^{(i')}(E).
\label{polaE}
\end{align}
From Proposition \ref{thm:conj3} we have the following formula, which is expressed as a hyperelliptic integral of second kind:
\begin{align}
& \Lambda (x+2\omega _k,E)=(-1)^{q_k}\Lambda (x,E) \cdot \exp \left( -\frac{1}{2} \int_{E_0}^{E}\frac{ -2\eta _k a(\tilde{E}) +2\omega _k c(\tilde{E}) }{\sqrt{-Q(\tilde{E})}} d\tilde{E}\right) ,
\label{hypellint}
\end{align}
for $k=1,3$, where $\eta _k =\zeta (\omega _k)$ and $\zeta (x)$ is the Weierstrass zeta function.
\begin{prop} (Hermite-Krichever Ansatz (\cite{Tak4} for Heun's equation)) 
Set
\begin{align}
& \Psi (x)=\prod _{i'=1}^M (\wp (x) -\wp (\delta _{i'})),\quad  \Phi _i(x,\alpha )= \frac{\sigma (x+\omega _i -\alpha ) }{ \sigma (x+\omega _i )} \exp (\zeta( \alpha )x),
\end{align}
for $i=0,1,2,3$, where $\sigma (x)$ is the Weierstrass sigma function. There exist polynomials $P_1(E),$ $\dots ,P_6 (E)$ such that, if $P_2(E') \neq 0$, then the function $\Lambda (x,E')$ is written as
\begin{align}
& \Lambda (x,E') = \frac{\exp \left( \kappa x \right) }{\Psi (x) } \left( \sum _{i=0}^3 \sum_{j=0}^{\tilde{l}_i-1} \tilde{b} ^{(i)}_j \left( \frac{d}{dx} \right) ^{j} \Phi _i(x, \alpha ) \right),
\label{Lalpha}
\end{align}
where the values $\alpha $ and $\kappa $ are expressed as
\begin{align}
& \wp (\alpha ) =\frac{P_1 (E')}{P_2 (E')}, \quad \wp ' (\alpha ) =\frac{P_3 (E')}{P_4 (E')} \sqrt{-Q(E')} , \quad \kappa  =\frac{P_5 (E')}{P_6 (E')} \sqrt{-Q(E')}.
\label{P1P6}
\end{align}
\end{prop}
Note that we have the periodicities
\begin{align}
& \Lambda (x+2\omega _k,E) = \exp(-2\eta _k \alpha +2\omega _k \zeta (\alpha ) +2 \kappa \omega _k )\Lambda (x,E) ,
\label{eq:monodHK}
\end{align}
for $k=1,3$.

We can obtain hyperelliptic-elliptic reduction formulae by comparing two expressions (Eq.(\ref{hypellint}) and Eq.(\ref{eq:monodHK})) of the monodromy.
\begin{prop} (Hyperelliptic-to-ellitic reduction integral formula (\cite{Tak4} for Heun's equation))
Set $\xi =P_1(E)/P_2(E)$. \\
(i) We have
\begin{align}
& \int _{\infty} ^{\xi } \frac{d \tilde{\xi }}{\sqrt{4 \tilde{\xi } ^3-g_2 \tilde{\xi } -g_3}} = -\frac{1}{2} \int _{\infty}^{E} \frac{a(\tilde{E})}{\sqrt{-Q(\tilde{E})}}d\tilde{E}.
\label{alpint}
\end{align}
(ii) Let $\alpha _0$ denote the value of $\alpha $ at $E=E_0$, where $E_0$ is the value satisfying $Q(E_0)=0$.
If $\alpha _0 \not\equiv 0$ (mod $2\omega _1 \Zint \oplus 2\omega _3 \Zint$), then $\kappa=P_5 (E)\sqrt{-Q(E)}/P_6 (E)$ is also expressed as

\begin{align}
 \kappa  & = -\frac{1}{2} \int  _{E_0}^{E} \frac{c(\tilde{E})}{\sqrt{-Q(\tilde{E})}}d\tilde{E} + \int _{\wp (\alpha _0)} ^{\xi } \frac{ \tilde{\xi } d \tilde{\xi }}{\sqrt{4\tilde{\xi }^3-g_2 \tilde{\xi }-g_3}} .
\label{kap123} 
\end{align}
\end{prop}

\begin{exa} We consider tha case $M=0$, $l_0=2$, $l_1=l_2=l_3=0$
Recall that the functions $\Xi(x,E) $ and $Q(E)$ are calculated as
\begin{align}
& \Xi (x,E)=E^2 +3\wp (x) E+ 9(\wp (x)^2 -g_2/4), \\
& Q(E)=(E^2 -3g_2)(E^3-9g_2E/4-27g_3/4).
\end{align}
The integral representation for solutions are expressed by substituting them into Eq.(\ref{eqn:Lam}).
The monodromy formula in terms of hyperelliptic integral are written as
\begin{align}
& \Lambda (x+2\omega _k,E)=\Lambda (x,E) \exp \left( -\frac{1}{2} \int_{\sqrt{3g_2}}^{E}\frac{ -6\eta _k \tilde{E} +\omega _k (2\tilde{E}^2-3g_2) }{\sqrt{-Q(\tilde{E})}} d\tilde{E}\right) ,
\end{align}
for $k=1,3$.
The function $\Lambda (x, E)$ admits an expression of the Hermite-Krichever Ansatz
\begin{align}
& \Lambda (x,E) = \exp \left( \kappa x \right) \left( \tilde{b}_0 \Phi _0 (x, \alpha ) + \tilde{b}_1  \frac{d}{dx} \Phi _0 (x, \alpha ) \right),
\end{align}
and the values $\alpha $ and $\kappa $ satisfy
\begin{align}
& \wp (\alpha )=-\frac{E^3-27g_3}{9(E^2-3g_2)} , \quad \kappa = \frac{2}{3}\sqrt{\frac{-(E^3-9g_2E/4-27g_3/4)}{(E^2-3g_2)}}.
\label{eq:alphalkappa2000}
\end{align}
The hyperelliptic-to-elliptic reduction integral formula for this case is written as
\begin{align}
& -\frac{1}{2} \int _{\infty}^{E} \frac{3\tilde{E}d\tilde{E} }{\sqrt{-(\tilde{E}^2 -3g_2) (\tilde{E}^3-\frac{9}{4}g_2\tilde{E}-\frac{27}{4}g_3) }} = \int _{\infty} ^{\xi } \frac{d\tilde{\xi } }{\sqrt{4\tilde{\xi } ^3-g_2\tilde{\xi } -g_3}}, \\
& \frac{1}{2} \int  _{3e_1}^{E} \frac{(\tilde{E} ^2-\frac{3}{2} g_2)d\tilde{E} }{\sqrt{-(\tilde{E}^2 -3g_2)(\tilde{E}^3-\frac{9}{4}g_2\tilde{E}-\frac{27}{4}g_3)}} = -\kappa + \int _{e_1 } ^{\xi } \frac{\tilde{\xi } d\tilde{\xi } }{\sqrt{4\tilde{\xi } ^3-g_2\tilde{\xi } -g_3}} ,
\end{align}
where $\xi =-(E^3-27g_3)/(9(E^2-3g_2))$ and $\kappa $ is defined as Eq.(\ref{eq:alphalkappa2000}).
These formulae reduce hyperelliptic integrals of genus two to elliptic integrals.
\end{exa}

\section{Finite-gap potential and WKB analysis} \label{sec:WKB}
The WKB analysis (the WKB approximation) appears in a semiclassical calculation in quantum mechanics,
and the WKB analysis is applied for the asymptotic analysis.
On the WKB analysis of the Schr\"odinger equation, a solution to the Schr\"odinger equation is obtained as a formal power series in $\eta ^{-1}$ by introducing a large parameter $\eta $, and it may have convergent expression by considering the Borel transformation for some cases (see \cite{AKT}).
We now consider a solution to the equation
\begin{equation}
\left(-\frac{d^2}{dx^2} + \eta^2 (v (x)-E) \right)\psi (x)=0,
\end{equation}
with a large parameter $\eta $.
We will find solutions in the form $\psi (x)=\exp(\int S(x,\eta )dx)$, 
\begin{equation}
S(x,\eta) = S_{-1} (x) \eta +S_0(x) +S_1(x)\eta^{-1}+S_2(x)\eta^{-2}+\dots .
\end{equation}
We set $Q=v (x)-E$. Then the function $S(x,\eta)$ would satisfy
\begin{equation}
S(x,\eta)^2+\frac{\partial S(x,\eta)}{\partial x}= \eta ^2 Q ,
\label{eq:S2}
\end{equation}
and we have the recursion formula
\begin{equation}
S_{-1}(x)^2=Q, \quad 2S_{-1}(x)S_j(x)= -\left( \sum_{k+l=j-1, k,l\geq 0} S_k(x)S_l(x) +\frac{dS_{j-1}(x)}{dx}\right) .
\end{equation}
We set $S_{odd}=\sum _{j\geq 0} S_{2j-1}(x) \eta ^{1-2j}$ and $S_{even}=\sum _{j\geq 0} S_{2j}(x) \eta ^{-2j}$.
Then we have 
\begin{equation}
S_{even}=-\frac{1}{2S_{odd} }\frac{\partial S_{odd} }{\partial x} , \quad \psi _{\pm } = (S_{odd})^{-1/2} \exp \left( \pm \int S _{odd} dx \right).
\label{eq:psi}
\end{equation}
Note that $\psi _{\pm }$ are formal solutions as power series in $\eta ^{-1}$.
Here we point out that Eq.(\ref{eq:psi}) resembles to an integral representation of solutions (Eq.(\ref{eqn:Lam})) in the theory of finite-gap potential.
In fact Eq.(\ref{eqn:Lam}) is written as Eq.(\ref{eq:psi}) up to a scalar multiplication by setting $ \tilde{S}_{odd} = \sqrt{-Q(E)}/\Xi (x,E)$.
Set $\tilde{S} = \tilde{S}_{odd} -(\partial \tilde{S}_{odd} /\partial x)/(2\tilde{S}_{odd} )$. 
Then an equality as Eq.(\ref{eq:S2}) follows from Eq.(\ref{const}).

We consider the case $v(x) =\wp (x)$, which is the case of Lam\'e's equation,
\begin{equation}
\left(-\frac{d^2}{dx^2} + \eta^2 (\wp (x) -E) \right)\psi (x)=0.
\label{eq:WKBLame}
\end{equation}
We have two methods for inventigating Lam\'e's equation, the finite-gap integration and the WKB analysis.
If $\eta ^2=l(l+1)$ for $l \in \Zint $, then the potential of Eq.(\ref{eq:WKBLame}) is algebro-geometric finite-gap and we can apply the finite-gap integration for Lam\'e's equation.
On the limit $\eta \rightarrow \infty$, we may apply the WKB analysis.
Fedoryuk \cite{Fed1,Fed2} investigated Lam\'e's equation, Lam\'e wave equation and Heun's equation by the WKB analysis and obtained error estimates for WKB approximation. 
On the WKB expansion of the solutions to Lam\'e's equation, we have
\begin{align}
& S_{-1}(x) = (\wp (x)-E)^{1/2} , \quad S_0 (x) = \frac{-\wp '(x)}{4(\wp (x)-E)} ,\\
& S_{1} (x) = -\frac{5\wp '(x) ^2}{32(\wp (x) -E)^{5/2}} +\frac{\wp''(x)}{8 (\wp (x) -E)^{3/2}}, \dots . \nonumber
\end{align}
We consider the leading term of $\eta ^{-1}$, i.e. $S_{-1}(x) $.
Then the monodromy with respect to the shift $x \rightarrow x +2\omega _j$ is written as
\begin{align}
& S(x+2\omega _j) = \exp \left( \pm \eta \int _{x}  ^{x + 2\omega _j} ( \wp (x)-E)^{1/2} dx +\eta^{-1}(\dots ) +\dots \right) S(x). 
\end{align}
We will sketch the distribution of the eigenvalues of Lam\'e polynomial for $\eta \rightarrow \infty$ and $l \in \Zint_{>0}$ ($\eta ^2=l(l+1)$).
Lam\'e polynomial is characterized as a non-zero doubly-periodic solution (i.e. $f(x+2\omega _i)/ f(x) \in \{ \pm 1\}$ for $i=1,3$) to Lam\'e's equation (see \cite{WW,Tak1}),
and it is essentially a polynomial in $z$ by setting $z=\wp (x)$.
If $l$ is a positive integer, then the number of the eigenvalues of the Lam\'e polynomial is $2l+1$ (see \cite{WW,BS}),
and the eigenvalue of the Lam\'e polynomial satisfies $\eta ^2e_1 \leq \eta ^2 E \leq \eta ^2 e_3$, (i.e. $e_1 \leq  E \leq e_3$) for the case $\omega _1 \in \Rea _{>0}$ and $\omega _3 \in \sqrt{-1} \Rea _{>0}$.
If $\eta ^2= l(l+1)$, then the eigenvalue $E$ of the Lam\'e polynomial satisfies
\begin{align}
& \eta \int _{\omega _3} ^{\omega _3 + 2\omega _1} (\wp (x ) -E )^{1/2} dx \in  \pi \sqrt{-1} \Zint , \quad (e_2<E<e_1),\\
& \eta \int _{\omega _1} ^{\omega _1 + 2\omega _3} (\wp (x ) -E )^{1/2} dx \in  \pi \sqrt{-1} \Zint , \quad (e_3<E<e_2),
\end{align}
by the asymptotics as $\eta ^{-1} \rightarrow 0$,
which follows from the periodicity of the monodromy.
It is shown under the condition $\eta ^2e_1 \leq \eta ^2 E \leq \eta ^2 e_3$, $\omega _1 \in \Rea _{>0}$ and $\omega _3 \in \sqrt{-1} \Rea _{>0}$ that a simply-periodic solution is also doubly-periodic.
It seems that $l$ pairs of eigenvalues merge by the WKB approximation.

We change a varible in the integral by setting $z=\wp (x)$.
The $m$-th eigenvalue $E \in (e_2, e_1)$ of the Lam\'e polynomial from the top should satisfy
\begin{align}
& \int _{e_3} ^{e_2} \sqrt{\frac{E-z}{(e_1-z)(e_2-z)(z-e_3)}} dz = \pi - \frac{m\pi }{\eta },
\end{align}
and the $m$-th eigenvalue $E\in (e_3 ,e_2)$ of the Lam\'e polynomial from the bottom should satisfy
\begin{align}
& \int _{e_2} ^{e_1} \sqrt{\frac{z-E}{(e_1-z)(z-e_2)(z-e_3)}} dz =  \pi - \frac{m\pi }{\eta }.
\end{align}
as $\eta ^{-1} \rightarrow 0$.
Hence the number of eigenvalues of the Lam\'e polynomial which is less than $E$ is  given by
\begin{align}
n(E)= \left\{ 
\begin{array}{ll}
\frac{\eta }{\pi } \left( \pi - \int _{e_2} ^{e_1} \sqrt{\frac{z-E}{(e_1-z)(z-e_2)(z-e_3)}} dz \right) ,& (e_3<E<e_2). \\
\frac{\eta }{\pi } \left( \int _{e_3} ^{e_2} \sqrt{\frac{E-z}{(e_1-z)(e_2- z)(z-e_3)}} dz\right)  , & (e_2<E<e_1), 
\end{array}
\right.
\end{align}
as $\eta ^{-1} \rightarrow 0$. Note that we used the formulae
\begin{align}
& \int _{e_3} ^{e_2} \frac{1}{\sqrt{(e_2-z)(z-e_3)}} =\int _{e_2} ^{e_1} \frac{1}{\sqrt{(e_1-z)(z-e_2)}} dz = \pi \\
& \left(\pi - \int _{e_3} ^{e_2} \frac{1}{\sqrt{(e_1-z)(z-e_3)}} dz\right)-\left(\int _{e_2} ^{e_1} \frac{1}{\sqrt{(e_1-z)(z-e_3)}} dz -\pi\right)=\pi . \nonumber
\end{align}
The density of the eigenvalues is written as
\begin{align}
\frac{1}{\eta } \frac{\partial }{\partial E} n(E) = 
\left\{ 
\begin{array}{ll}
\frac{1}{2\pi } \left( \int _{e_2} ^{e_1} \frac{1}{\sqrt{(e_1-z)(z-e_2)(z-e_3)(z-E)}} dz \right) ,& (e_3<E<e_2), \\
\frac{1}{2\pi } \left( \int _{e_3} ^{e_2} \frac{1}{\sqrt{(e_1-z)(e_2- z)(z-e_3)(E-z)}} dz\right)  , & (e_2<E<e_1). 
\end{array}
\right.
\end{align}
Hence we recover the formula given by Borcea and Shapiro \cite{BS} by the completely different method.
We decompose the integral as
\begin{align}
& \frac{1}{2\pi} \int _{e_3} ^{e_2} \frac{dz}{\sqrt{(e_1-z)(e_2-z)(z-e_3)(E-z)}} \\
& = \frac{1}{2\pi\sqrt{(e_1-E)(E-e_3)}} \int _{e_3} ^{e_2} \frac{dz}{\sqrt{(e_2-z)(E-z)}} \nonumber \\
& + \frac{1}{2\pi}\int _{e_3} ^{e_2} \frac{dz}{\sqrt{(e_2-z)(E-z)}} \left( \frac{1}{\sqrt{(e_1-z)(z-e_3)}}-\frac{1}{\sqrt{(e_1-E)(E-e_3)}} \right) . \nonumber
\end{align}
Then the asymptotic of the density of the eigenvalues as $E \rightarrow e_2+0$ is written as
\begin{align}
& \frac{1}{2\pi } \left( \frac{\log \left( \frac{4(e_2-e_3)}{(E-e_2)} \right)}{\sqrt{(e_1-e_2)(e_2-e_3)}} + \frac{\log \left( \frac{4(e_1-e_2)}{(e_1-e_3)} \right)}{\sqrt{(e_1-e_2)(e_2-e_3)}} \right) =\frac{1}{2\pi } \frac{\log \left( \frac{16(e_1-e_2)(e_2-e_3)}{(e_1-e_3)(E-e_2)} \right)}{\sqrt{(e_1-e_2)(e_2-e_3)}}.
\end{align}
It is also shown that the asymptotic as $E \rightarrow e_2-0$ is written as
\begin{equation}
\frac{1}{2\pi }\frac{\log \left( \frac{16(e_1-e_2)(e_2-e_3)}{(e_1-e_3)(e_2-E)} \right)}{\sqrt{(e_1-e_2)(e_2-e_3)}}.
\end{equation}
Hence the density of the eigenvalues as $E\rightarrow e_2$ has logamithtic singularity.
If $e_2=e_3$, then the density is written as $1/(2\sqrt{(e_2-e_1)(E-e_2)})$ and the type of the singularity at $E=e_2$ changes.

We have another way to insert a large parameter.
we regard the eigenvalue $E$ as a large parameter while fixing the potential.
Namely, we consider the Schr\"odinger equation
\begin{equation}
\left(-\frac{d^2}{dx^2} +  (v (x)-E) \right)\psi (x)=0,
\end{equation}
by setting $E=-\eta ^2$,
\begin{equation}
\psi (x) = \exp ( \eta x  +\psi _1(x)\eta^{-1} +\psi _2(x)\eta^{-2}+\psi _3(x)\eta^{-3} +\psi _4(x)\eta^{-4}+ \dots ).
\end{equation}
and we regard $\eta $ as a large parameter while fixing the potential.
We now restrict to the case of Lam\'e's equation $v(x)=l(l+1)\wp (x)$.
Then the coefficients $\psi _1(x), \psi _2(x) ,\dots $ are expressed as
\begin{align}
& \psi _1(x)= -\frac{l(l+1)}{2} \zeta (x), \quad \psi _2(x) = -\frac{l(l+1)}{4} \wp (x) , \\
& \psi _3 (x) = -\frac{l^2(l+1)^2}{96}g_2 x+\left(-\frac{l^2(l+1)^2}{48}+\frac{l(l+1)}{8}\right) \wp '(x) , \nonumber \\ 
& \psi _4 (x) = \frac{l^2(l+1)^2}{96}g_2 +\left(\frac{l^2(l+1)^2}{48}-\frac{l(l+1)}{16}\right) \wp ''(x) , \dots . \nonumber
\end{align}
Hence we have the asymptotics of the monodromy
\begin{align}
\psi (x+2\omega _i)= \psi (x) \exp \left( 2\omega _i \eta - l(l+1) \eta _i  \eta^{-1} -\frac{l^2(l+1)^2 g_2 \omega _i}{48}\eta^{-3} +\dots \right) ,
\label{eq:psimonodGV}
\end{align}
where $\eta _i= \zeta (\omega _i)$ $(i=1,3)$.
Grosset and Veselov \cite{GV1} studied the asymptotics of the monodromy by applying the elliptic Faulhaber polynomials, which is related to the soliton theory of the KdV equation.
They considered the asymptotics of the densities of states for  Lam\'e's equation,
and applied it for computation of the polynomials $-2\eta _k a(E) +2\omega _k c(E) $ in Eq.(\ref{hypellint}).
Eq.(\ref{eq:psimonodGV}) is also obtained as a consequence of the paper \cite{GV1}.

We expect further results among the finite-gap potential, WKB analysis and Heun's (Lam\'e's) equation.

\subsection*{Acknowledgements}

The author would like to thank Professor Takashi Aoki for providing me an opportunity to give a talk on the workshop "Algebraic Analysis and the Exact WKB Analysis for Systems of Differential Equations" (December 2006, RIMS, Kyoto).
Thanks are also due to Professors Naoto Kumano-go and Susumu Yamazaki.
He thanks Professor Boris Shapiro and Dr. Tatsuya Koike for discussions.


\end{document}